# Data Mining in Electronic Commerce

**David L. Banks and Yasmin H. Said**

*Abstract.* Modern business is rushing toward e-commerce. If the transition is done properly, it enables better management, new services, lower transaction costs and better customer relations. Success depends on skilled information technologists, among whom are statisticians. This paper focuses on some of the contributions that statisticians are making to help change the business world, especially through the development and application of data mining methods. This is a very large area, and the topics we cover are chosen to avoid overlap with other papers in this special issue, as well as to respect the limitations of our expertise. Inevitably, electronic commerce has raised and is raising fresh research problems in a very wide range of statistical areas, and we try to emphasize those challenges.

## 1. INTRODUCTION

Electronic commerce is changing the face of business. It allows better customer management, new strategies for marketing, an expanded range of products and more efficient operations. A key enabler of this change is the widespread use of increasingly sophisticated data mining tools.

The Department of Commerce commissioned a study of 2003 economic data (U.S. Census Bureau, 2005). It showed that e-commerce, on a percentage basis, outperformed all four major economic sectors in 2002–2003. For the manufacturing sector, 21.2% ($843 billion) of the activity (as measured in total sales dollars) was classified as e-commerce. For the merchant wholesalers sector, e-commerce sales were 16.9% ($730 billion) of total sales; for retail trade it was 1.7% ($56 billion) and for selected service industries it was 1% ($50 billion). The dominant component was business-to-business activity. These trends can only have increased since 2003.

Nonetheless, reviews and research in this area are handicapped by the proprietary nature of the data and the algorithms. A great deal of effort is being expended in this area, but most of it is secret. Certainly Amazon, Google and Microsoft are deeply engaged in statistical research, and in time the broader research community may learn more about their findings, but for now, all this paper can really attempt is to lay out the main strategies in the relevant areas.

*David L. Banks is Professor, Institute of Statistics and Decision Sciences, Duke University, Durham, North Carolina 27708, USA e-mail: banks@stat.duke.edu. Yasmin H. Said is Lecturer, Department of Applied Mathematics and Statistics, Johns Hopkins University, Baltimore, Maryland 21218, USA e-mail: ysaid99@hotmail.com.*

In that context, we try to survey the contributions that data mining makes to e-commerce, and point out legal, social and commercial issues that arise in its use. We also address research issues that arise in such applications. After a short introduction to the development of e-commerce, we treat the paradigm example of web browsing and search engines in Section 2, and customer relationship management (CRM) in Section 3. Section 4 concludes with a discussion of business strategies for the future.

### 1.1 History of E-Commerce

It seems presumptuous to discuss the history of a technology change that is barely a decade old. However, it is helpful to review this history to emphasize the rate and magnitude of the e-commerce transition. In addition, it is important to point out how unforeseen consequences of the technology change are driving the research challenges in today's applications.

Until about 1994, electronic commerce was not web-based. The term referred to the use of computers and telecommunications to automatically forward and process commercial documents, such as invoices and inventory requests. A major figure in this area was Ross Perot, who founded the Electronic Data Systems Company (EDS) in 1962 (the initial goal was to streamline municipal parking ticket billing). The core technology was an industry standard (ANSI X12, X.400) for electronic data interchange (EDI), which allowed communication between computers, originally through shipment of magnetic tapes. Until the 1990s, early versions of electronic commerce focused on business-to-business (B2B) transactions, because personal computers were relatively rare and the EDI systems were expensive.

However, the infrastructure was growing rapidly—Figure 1 shows the estimated growth in both personal computers and the size of the Internet. Some large businesses foresaw the potential (and some, notably Microsoft, almost came to the table too late), but much of the growth was driven by amateurs. One of the remarkable qualitative differences in the Internet revolution, as compared to previous innovations, was that it lowered the capitalization threshold for new businesses. Netscape, Yahoo!, Google, Travelocity, Mapquest, Amazon and e-Bay all began on a very small scale, and their primary assets were the technology skills of the early entrepreneurs.

At first, most e-commerce start-ups had simple business models rooted in traditional bricks-and-mortar perspectives. A typical paradigm was to have customers place orders over the Internet, or to sell advertisement space (pop-ups and banners) during web browsing sessions. However, emergent complexity soon began to dominate the evolution, and business biodiversity blossomed.

In the case of Amazon and similar companies (such as Travelocity and Expedia), the initial plan was to reduce costs by cutting out storefronts. The merchandising was modeled after the old Sears and Roebuck catalog sales, except that people would place orders over the Internet and receive deliveries by Federal Express and/or conventional post. Similarly, the search engine companies (and MapQuest) initially saw themselves as high-tech telephone books, with revenue generated from advertisements, but in both cases new business opportunities arose as the result of applying data mining algorithms to recommend books or hotels to customers, or to target advertising more effectively, based on activity histories.

After this first wave of e-commerce business strategies, secondary growth in niche markets began to boom. These ranged from on-line dating services to medical advice (drkoop.com) to character assassination (RateMyProfessors.com). Computer security software became a necessary enabler, as well as spam filters and technology such as PayPal to facilitate small-scale commercial purchases. Gray commerce became pervasive, notably in pornography and services such as Napster that support music file-sharing. The norms and laws of society are struggling to adapt to new circumstances; this is complicated by the fact that businesses can easily be located outside of national jurisdiction where different rules apply. Some nations are more advanced than others with respect to specific aspects of e-commerce; for example, the United Kingdom is moving quickly to-ward universal Internet banking and Japan is pioneering a wide range of popular innovations.

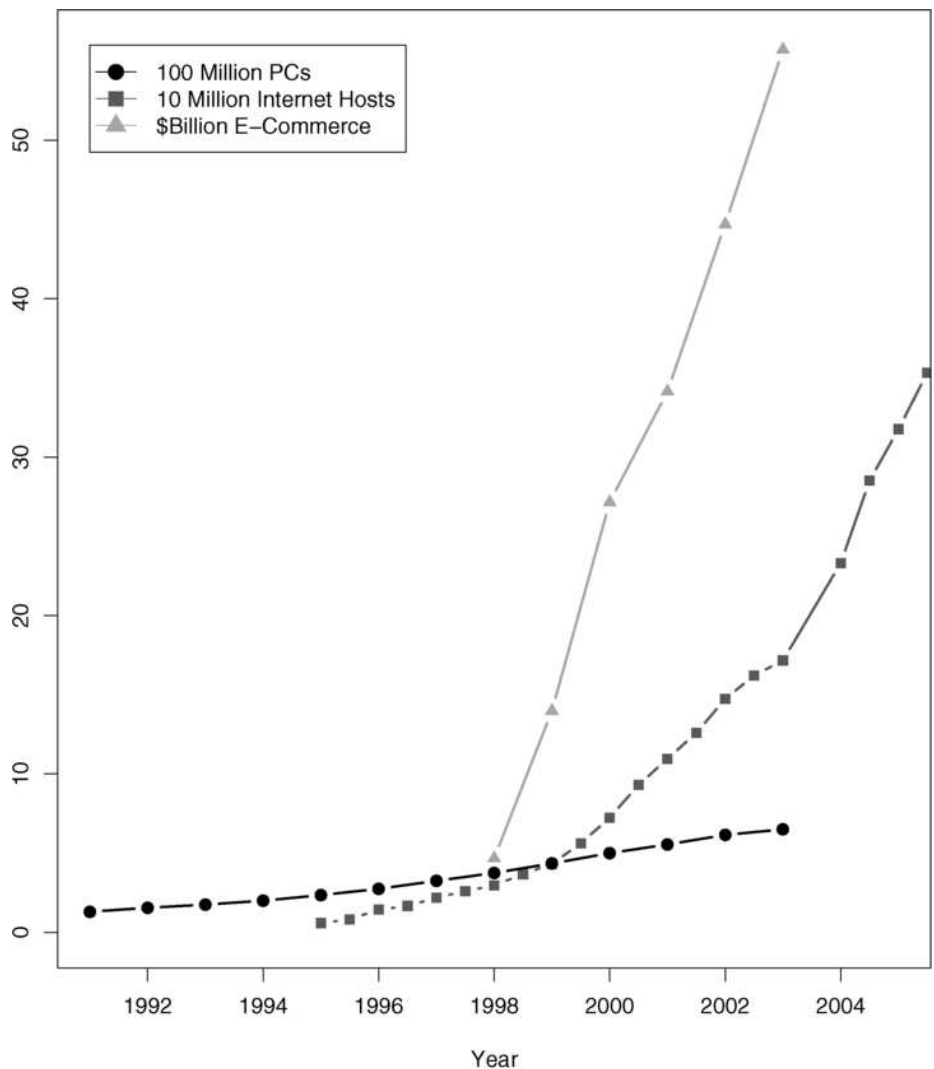

FIG. 1. *Figure* 1 *shows how growth in the estimated number of PC-owners preceded the explosive growth in e-commerce, which roughly tracks the growth in the estimated number of servers. The data on the number of PC owners is from the International Telecommunications Union (2004), the data on e-commerce is from* www.census.gov/eos/www/historical/2003ht.pdf, *and the data on servers is from the Internet Systems Consortium (*www.isc.org*).*

More recently, businesses that aggregate data for resale have appeared. These firms gather and organize data so that it can be analyzed by others (in contrast to data collected to meet specific business objectives). Lexis–Nexis–Seisint and ChoicePoint are two examples of firms that collect electronic data about people from multiple sources, analyze it and sell the material to others. Seisint and ChoicePoint initially worked primarily with credit institutions, who were concerned with identity authentication, credit security and the ability to distinguish risk and profit categories among potential clients. Various law enforcement agencies are also purchasing data about people, presumably to search for terrorists (among other things). This has redoubled concerns about privacy rights, highlighting issues that had been irrelevant when data were hard to collect and difficult to synthesize. Hui and Png (2006) examined the economics of privacy, and pointed out perverse incentives and asymmetric information.

Much of the growth and diversification synopsized above depends on data mining methods (of varying degrees of sophistication), but e-commerce also allows businesses to access new data streams that inform management in ways that were not previously possible. One example is the on-line auctions at e-Bay, which provide data on how much money customers are willing to pay for a product. Another example is clickstream analysis, which provides clues on how people weigh information in making purchase decisions. A third is the expanded ability to rapidly identify product defects (e.g., monitoring user groups quickly found the legacy code flaw in *Grand Theft Auto*: *San Andreas*). A fourth example is the application of social network analysis to fine-structure data on the relationships among people and businesses. Bapna, Goes, Gopal and Marsden (2006) described automated data collection methods in several areas of e-commerce that enable low-cost tests of economic theories. More strongly, Allen, Burk and Davis (2006) laid out ethical and practical guidelines to which academic researchers should adhere when collecting e-commerce data.

Before finishing this historical review, collegiality demands an acknowledgment that mathematical disciplines other than statistics have contributed to the phenomenal growth of e-commerce. In particular, progress would have been stymied without encryption methods that enabled secure transactions. Most of these methods are based on trap-door codes discovered much earlier: the RSA system (cf. Rivest, Shamir and Adleman, 1978) has security equivalent to the difficulty of factoring the product of two large primes; the BBS system (cf. Blum, Blum and Shub, 1986) is based on the complexity of modular exponentiation. Information theory and queuing theory have also made crucial contributions.

In addition to these key enablers, mathematics has led to new services. For example, MapQuest depends strongly on sophisticated algorithms for generating driving directions, which draw heavily on graph theory and optimization; various on-line games (e.g., *EverQuest*, a popular multi-user dungeon) require sophisticated mathematics for visualizations; and network theory is embedded in the recent emergence of facilitated professional connection-building, such as LinkedIn.com (who can guess where such click-streaming will go?).

### 1.2 Research Issues in E-Commerce

People want to know about how the new technology is changing business. Hot topics are spam, social network analysis, text mining of discussion boards, pricing and fraud detection. However, despite the new low cost of data and the urgent interest, there are major obstacles to research.

First, for nearly all U.S. academic researchers, an institutional review board must approve human subject studies. This requirement grew from concerns that arose in several famous experiments in psychology and sociology, and also from the need to ensure fair practice in the conduct of clinical trials (cf. *The Belmont Report*, National Commission for the Protection of Human Subjects, 1979). The relevance of these issues to e-commerce was confirmed by a notorious study undertaken at Carnegie Mellon; Martin Rimm, an undergraduate, worked with his advisor Marvin Sirbu to analyze patterns in Internet pornography use by the faculty, staff and students, as well as data from commercial pornographic web sites (cf. Rimm, 1995).

A second obstacle is the proprietary nature of most e-commerce data. The data are an asset; the companies want to mine it themselves and to keep it from their competitors. There are some exceptions: e-Bay bid sequence data are available for 15 days after the close of the auction, and they have been analyzed by Shmueli and Jank (2006); Amazon price data are freely available, as analyzed by Ghose and Sundararajan (2006). Even when open data exist, they can be hard for an outsider to know the timings of promotional campaigns and web site changes that affect consumer behavior. Also, the tools and assumptions that are used for mining are often tuned to the problem (especially in applications such as text mining). In principle, there may be a way forward for some kinds of data sharing: statisticians have developed methods to analyze multiple data bases without combining them directly (cf. Karr, Lin, Sanil and Reiter, 2005). This allows companies to learn faster from their joint data without having to actually disclose any of their raw information. However, the methods for secure analysis of distributed data bases apply to simple summary statistics and multiple linear regression: these low-level analyses will not satisfy either the wishes of academic researchers or the needs of business managers.

The third obstacle to e-commerce research is almost trite. For many applications, the quality of the data is poor and enormous preprocessing is required. There are few statistical fixes for this. Some approaches include automatic reconciliation, which is based on record linkage strategies implicit in Fellegi and Sunter (1969), use of highly robust analyses whose breakdown point can cope with large amounts of bad data (cf. Donoho and Huber, 1983) or visualization methods for cleaning data (Karr, Sanil and Banks, 2006). However, the cost of retroactively cleaning data sets is usually not commensurate with the value obtained from their analysis. It is better to focus on improving future data quality than to comb out the tangles in the past.

Nonetheless, data mining is the future in e-commerce. Kohavi, Mason, Parekh and Zheng (2004) laid out a set of research issues related to retail e-commerce, a satellite meeting on web mining and web usage analysis was part of the 2004 Annual Conference on Knowledge Discovery and Data Mining, and Darryl Pregibon, now at Google, gave a recent talk at the University of Maryland in which he indicated that Google was using Bayesian methods for spelling correction, Bayesian nets for clustering queries and log-linear models for translation among languages. These applications, and those discussed in more detail in this paper, illustrate the first fruits of modern statistics in this arena.

## 2. WEB BROWSING

The breakthrough that made e-commerce possible is the web browser. Invented in 1990 by Tim Berners-Lee, a CERN physicist, to support his other invention, the linked document set that became known as the World Wide Web, it quickly became the paradigm killer application. The web browser displays and manipulates documents hosted by web servers or held in a file system, using Internet protocols to navigate among links. It is most useful when linked to a search engine, allowing users to find documents with specified key words. In particular, it allows customers to find sites that will sell them the goods they want.

Search engines do not actually search the whole World Wide Web. Instead, they search a data base of preselected sites that have been preprocessed to create an index that summarizes content features. This data base gets updated every few hours or days, with information obtained by spiders (or 'bots or intelligent agents). These spiders move from web site to web site, deciding which ones to select and index. When a user submits a query, the data base is searched and the sites whose indexes show a close match are returned as hits, ranked according to some criterion for relevance.

One of the ways in which commercial search engines compete is in terms of the size of their data base, but these numbers change frequently, are difficult to verify, may include duplicates and are probably not very relevant to overall performance. Nonetheless, to give a sense of the scale, the numbers in Table 1 were reported as being current on November 11, 2004 by a specialist search engine web site (Sullivan, 2004). Note that Google started the practice of counting "partially indexed" web pages and most other companies have followed. These are pages that are only known through links and have not actually been summarized for searching.

Page depth may be a more important aspect of performance than size; it indicates how much of the web page is used for indexing. Google's rule is that if a page exceeds 101K of text, then Google uses only the first 101K, ignoring everything else. Competing search engines claim to index more (and this is testable: one can create a large web page, bury the key word "aardvarkpants" at a certain depth and see whether it gets indexed).

Search engine performance is subtle. Different engines take different approaches (Google can be quite unstatistical: it sells priority placement on retrieval lists, which may overstate the relevance, but appears to work well for commercial applications). Bradlow and Schmittlein (2000), in a now somewhat dated article, described issues in performance measurement, and performance is closely related to the size of the Web—an important estimation problem discussed by Dobra and Fienberg (2003).

In the context of web browsing, two main data mining issues arise. The first is to find algorithms for text retrieval: the search engine with the competitive edge is the one that is best able to find relevant matches to a human inquiry (which is often ambiguous and casually worded); this entails record linkage and cluster analysis. The second issue is the development of tailored marketing, so that pop-up ads (the main revenue stream) are targeted to the customers most likely to respond.

TABLE 1
*Size and depth for some popular engines*

| Search engine | Reported size | Page depth |
| --- | --- | --- |
| Google | 8.1 billion | 101K |
| MSN | 5.0 billion | 150K |
| Yahoo! (est.) | 4.2 billion | 500K |
| AskJeeves | 2.5 billion | 101K+ |

### 2.1 Text Retrieval

Analyzing the interplay of factors that affect the performance of information retrieval systems is a hard problem. Early search engines had users type in key words, and then applied stemming and a bag-of-words model to assess relevance. Modern search engines have expanded capability: AskJeeves users express their information need in natural language; some search engines do cluster analyses of documents that seem relevant, so that the user can iteratively eliminate whole groups of impertinent documents.

2.1.1 *Stemming.* "Stemming" refers to the preprocessing step in which a key word (or a word in the retrieved document) is reduced to a generic form. If a user typed in "neural networks," he or she probably would not want to exclude documents that mention only "neural nets" or "neural network." The stemming process reduces the morphological variants of a word (plurals, tenses, gerunds, prefixes and so forth) to a common form. This common form need not be a proper word—"multivariate" and "multivariable" might both stem to "multivaria." Stemming generally improves the quality of the search and has the additional benefit of reducing the number of terms in the dictionary (usually by a factor of 5 or so).

There are obvious dangers in stemming. "Hamming" is not usefully stemmed to "ham," although "spamming" should be reduced to "spam." Homographs are an issue ("fleet of foot" versus "fleet of ships"). In addition, different languages require different stemmers, as do specialized subjects with nomenclature conventions (e.g., chemistry). To fully embrace the stemming strategy, one should include synonyms and equivalent phrases, but this quickly becomes labor-intensive, and overstemming can make it difficult to search for exact quotes or to check student essays for web plagiarism.

The stemming community has invented a number of open-source algorithms (the Paice/Husk algorithm, the Porter algorithm, etc.), many of which are currently available at www.comp.lancs.ac.uk/computing/research/stemming/general/index.htm. Necessarily, there are also many procedures for comparing these algorithms (cf. Paice, 1996); these include modified Hamming distances, dictionary compression factors and measures of correlation on sets of challenging words. The general view is that understemming is preferable to overstemming for search engines (in part because most users are able to add terms that refine the search when needed).

Commercial search engines have invested substantial resources in building stemming rules and they probably review user records of real searches to identify improvements. Data mining methods would be a natural part of this effort. Specifically, researchers would take some measure of stemming performance (say accuracy on a set of challenge words) and then use feature selection to build rule sets with good generalizability to most applications. In the same way that boosting iteratively upweights misclassified cases, boosted stemmers could tune their rule sets to handle the many exceptional cases that occur (and that freshly arise in living languages).

2.1.2 *The "bag-of-words": model and extensions.* One way to index a web page is to do it by hand, scoring according to prescribed set of categories. This approach is partially used by Yahoo! and it works well for popular topics, but it is slow and costly to update. Alternatively, one can automatically make a list of all the words in the document, perhaps weighted by frequency, and standardize by the vocabulary diversity (so as not always to retrieve the dictionary).

The latter approach is called a *bag-of-words* model because it ignores word order, titles and section headings, and many other clues about the significance of the term. From a data mining perspective, it corresponds to a naive Bayes classifier because it does not use the "correlation" structure, treating all (stemmed) words in the same way—any random permutation of the text would produce the same match. Bickel and Levina (2004) explained why naive Bayes classifiers can work so well (the cumulative inaccuracy from estimating the full covariance structure overwhelms the signal), and the same intuition applies here. Their suggestion to estimate only the most significant covariances corresponds to using only the most important contextual information. (Hand and Yu, 2001, gave other reasons why naive Bayes classification is successful and they also apply.)

Smarter systems try to use important context. Google puts extra weight on pages that have many in-links, since these appear to be popular. Also, most engines give more weight to search term matches in the title than in a footnote, and perhaps terms that occur early in the text count more. Additionally, sequences of within-sentence capitalization can signal a place or person's name: "Bill Gates" is treated differently from "We shall bill gates to his account." These methods use ad hoc rules, many of which are related to those used in natural language processing (e.g., AskJeeves; also, the Lycos search engine was based in part on FERRET, an earlier system that used word context; cf. Mauldin, 1991).

In contrast to this "ad hockery," a number of retrieval strategies are explicitly statistical. The Inktomo search engine formerly used logistic regression to predict the probability that a document was relevant as a function of key word matches (and probably other covariates, such as propinquity within sentences or use in the title). Excite's concept-based search engine uses latent semantic indexing, which is a singular value decomposition of the matrix whose columns correspond to terms and whose rows correspond to web pages (cf. Dumais, 1991). This is an attempt to reduce the dimension of the problem and it puts documents that have similar word choices near together in the reduced space (e.g., some web sites might refer to the "Mann–Whitney test," others to "Wilcoxon test" and others to the "Mann–Whitney–Wilcoxon test"; the third class provides a bridge that enables all three groups to be co-relevant, even though the first two have no terms in common).

There are many other methods discussed in the open literature. Hidden Markov models for generation of search terms apply data mining tools to address the problem (Miller, Leek and Schwartz, 1999). Schapire, Singer and Singhal (1998) used boosting methods to classify documents as relevant or not. Given the training of some of the researchers at Yahoo! and Google, it is likely that the closed literature employs tree-based rules and support vector machines for this classification problem.

2.1.3 *Comparing search engines.* For most e-commerce businesses the details of the search engine algorithms are proprietary; sometimes even their general strategy is opaque. Attempts to make head-to-head comparisons of these black boxes is difficult.

One basis for comparison is the coverage and overlap of the search engine's data bases of indexed pages. A classic and incomplete study on this was done by Lawrence and Giles (1999); regrettably, it is dated and not very statistical. A relatively new data mining technique should be useful here; it is based on the multiple systems methods developed by Ball (2003) in the context of human rights work (Ball's paper is based on joint research with Fritz Scheuren and Herbert Spirer). Their problem was to use data bases of rights violations (say disappearances) collected by different organizations (e.g., Human Rights Watch, the Catholic Church and local advocacy groups) to obtain an estimate of the number of violations that were never reported, based on the amounts of overlap among the data bases.

The intuition behind multiple systems estimation reflects both capture–recapture estimates of population size and the Good–Turing estimate (Good, 1953) for the number of unobserved species on an island (one knows the number of species observed exactly once, twice, etc.; a good nonparametric estimate of the number of unseen species is the number of species observed only once). The theory assumes independent sampling, which breaks down for both human rights violations and web pages, but a more complex model can give reasonable results. That model assumes that, conditional on the parameters, data are independent; for rights data, the parameters might indicate province and year, while for search engines, they might indicate content areas (business, scholarship, blogs). One especially nice feature of multiple systems estimates is that they can be done at the query level: one can look at the documents retrieved by each search engine for a specific query, count the overlaps, and get an estimate of the number of unseen documents and the coverage of each engine for that specific subject.

A better basis for comparison is performance. This is usually measured in terms of recall and average precision. Recall is the fraction of relevant documents on the Web that were found by the search engine, while precision assesses the fraction of a system's retrieved documents that are actually relevant; these ideas are related to Type I and Type II error or to sensitivity and specificity. Essentially, each system must balance the risk of excluding a relevant document against the risk of including an irrelevant one.

To estimate recall, one needs an estimate of the number of relevant documents on the Web. In principle, this can be found through multiple systems estimation, where each search engine generates its own list and the overlaps enable an estimate of the total. To estimate the average precision (AP), the formula is

$$\mathrm{AP} = R^{-1} \sum_{d \in D} \frac{\#\{\text{relevant documents retrieved at or before document } d\}}{\#\{\text{number of documents retrieved at or before document } d\}} \times I(d),$$

where $R$ is the number of relevant documents, $D$ is the set of all documents, $\#\{\cdot\}$ denotes the cardinality of the argument set and $I(d)$ is an indicator function which takes the value 1 if $d$ is a relevant document and is 0 otherwise. This measure was proposed by Harman (1994) and has been widely employed ever since. One can view this as the area under the step function that represents the cumulative proportion of relevant documents found in the ranked list (where documents are ranked from most relevant to least).

A main problem with both recall and average precision is the need to determine whether a document is relevant. Some users have richly associative minds; to them, lots of documents are relevant. Other users are more tightly focused.

The best current studies in the public domain are from the Text Retrieval Conference (TREC) bake-offs run by the National Institute of Standards and Technology. These began in 1992 and compare (mostly noncommercial) search engines on test-bed collections of news wire documents. The study protocol is to have retired CIA analysts craft queries (e.g., "Black Bear Attacks") and then go through each document in the collection to determine its relevance. These determinations become the gold standard against which many search engines (mostly written by university researchers) are compared. It is worth noting that for some documents, the analysts disagree among themselves about relevance.

From a data mining perspective, a number of methods have been tried. Banks, Over and Zhang (1999) used visualization, cluster analysis, Mandel's bundle-of-lines model, rank correlations and multidimensional scaling (MDS); aside from a few rather clear-cut queries, the large number of outliers made these methods unsatisfactory. Liggett and Buckley (2005) confirmed that MDS can work well for some queries. Finally, ensemble classification systems, where different systems vote on relevance or have their judgments weighted as in "stacking" (cf. Hastie, Tibshirani and Friedman, 2001, Chapter 8), provide strategies for improving performance. Many systems of this kind have been studied by TREC [and related programs, such as the ARPA (Advanced Research Projects Agency) program TIPSTER].

### 2.2 Tailored Marketing

Tailored marketing is a key revenue stream for commercial search engines, as well as for MapQuest, WhitePages, Dictionary.com and many other free online services. It also underlies recommender systems, such as the one used by Amazon.com to suggest books that may match a customer's interests. The intent is that the information that becomes available during an Internet session, or that is cumulatively available across sessions when cookies, purchase histories or other tracking systems are employed, enables advertisers to target offers toward consumers who are predisposed to buy. Someday, it may even forecast a person's price point, so that the price of the good is tuned to an individual's utility function (a modern e-souk!). Amazon is trying something like this already: to see this, check the price of several books, then delete the Amazon cookies and recheck the prices.

Some people are uncomfortable with such salesmanship; others see the technology as a potential boon that could someday reduce the overhead (in time and cost) of shopping. The issues touch closely on privacy concerns, as discussed by Fienberg (2006).

From a data mining standpoint, the primary issues are technical. One wants to build a model that uses available information about the Web user to intelligently guess the next ad to show. The structure of such problems is new and exciting. Dynamic programming is one component; it determines the best ad to offer next. Classification models are another component; these determine the probability of making a purchase as a function of all available information. Intelligent exploration is a third component; it is a good idea to offer a small number of ads at random to test the model and to learn more about the consumer. How these ingredients should be stitched together is an open research area (although it is likely that very smart people are working on this behind the fire walls of corporate research).

For dynamic programing formulation, this is a continuous time, discrete state problem with a potentially infinite horizon. The states include such things as seeing the ad, closing the ad, deciding to buy, seeing a new ad and terminating the session. State transitions occur at random times, but these intervals offer useful clues to the marketer: a long interval before closing a pop-up may indicate that the user is reading carefully and is tempted to buy, so a similar product should be offered in the next ad presented. The goal of the dynamic programing problem is to present a sequence of ads that enables the classification rule to learn quickly about the customer and thus to maximize the total purchase.

The user's clickstream is the raw data for the dynamic programming problem. The clickstream is a time-stamped sequence of the on-line activity, which can also help in web site design and estimation of price points. Moe and Fader (2004) described a statistical analysis of clickstream data regarding purchase decisions; Chatterjee, Hoffman and Novak (2003) use clickstream data to assess banner advertising. Some estimation methods for choice models in clickstream data were described by Sismeiro and Bucklin (2004).

Intelligent exploration is needed to learn about changing situations. Marketing evolves: new ads are developed, new pricing or financing options emerge and marketers want to assess the impact. Also, customers change: when someone has a child or moves to a new city, their shopping profile alters. In addition, consumer groups also change: new music fads arise, new vacation spots become hot and suburban commuters may suddenly develop a serious interest in fuel efficient cars. To keep up with these moving targets, smart marketing should frequently dangle new bait. This experimentation involves calculated risks, and one wants to ensure that the costs of losing a customer by trying a misconceived ad are outweighed by the value of the information obtained. Experiments should be designed to target specific customer segments to the extent that segment information is available.

Classification is a standard problem in data mining. In e-commerce applications, the appropriate method [random forests, SVMs (Support Vector Machines), logistic LASSO, etc.] depends sensitively on the web site, the type of ads that are available and the kinds of information that can be collected. For Amazon.com, the recommender system has access to previous purchases and books inspected. This is a richer source of information than Dictionary.com can access (they only know the word someone looks up in the current session, unless they have cookies). Some businesses obtain additional information from data warehouses, such as ChoicePoint, and this enables their experts to build classification rules that are highly individualized.

To be a little more concrete, suppose someone uses a web browser to access a site about, say, professional wrestling. It is trivial and obvious for the search engine to advertise event tickets and relevant tee shirts. If the marketing models are sophisticated, they might proffer country music CDs, pickup trucks and GED programs, and if the system knows that the user recently purchased *Have a Nice Day* (by Mick Foley a/k/a Mankind, Dude Love and Cactus Jack), then it might use that information to recommend *The Rock Says* (by Joe Layden a/k/a The Rock).

Such suggestions require a classification model built from market research data about the kinds of products that appeal to different clusters of the population. The dynamic programing component picks which ads to offer in which order, based on incoming information on what the user opens and closes, how long the user appears to study various ads, and the payoff from different kinds of purchases. To ensure that the models are tracking societal change and personal variation, the exploration component should sometimes, rarely and randomly, probe with an unorthodox recommendation; for example, the wrestling fan might see a pop-up for a statistics book.

## 3. CUSTOMER RELATIONSHIP MANAGEMENT

One of the new old things that e-commerce can do is customer management. When shops were small and communities were insular, each proprietor handled customer relations automatically. As mercantile empires formed, the customer/merchant relationship became impersonal. E-commerce now allows the possibility of recovering some of the individualized service that can cement return business.

Customer relationship management is a generalization of the tailored marketing discussed in Section 2.2. It goes beyond advertising to include all aspects of the customer experience: contact, billing, retention, help desks and even holiday e-cards. Successful use requires detailed files on each customer; one uses data mining to anticipate the kind of relationship that specific people want.

Early efforts at this were based on market segmentation. Businesses attempted to discover clusters of consumers who were similar, and then would develop payment plans, ad campaigns, special discounts and other policies designed for each cluster (especially the most profitable). The data mining tool used for this was cluster analysis, and the most famous commercial pioneer is Claritas, which used Census data to identify 64 "clusters" of consumers, with shorthand descriptors such as "kids and cul-de-sacs" or "money and brains" or "back country folks." Clusters get revised periodically to reflect important changes; for example, Claritas recently added the cluster "young digerati" to reflect the important technophile segment.

Customer relationship management can use such clusters to build models. One approach to dimension reduction is to build a separate model for each cluster; in this way, variables that are significant for the consumer behavior of "back country folks" but irrelevant to the "young digerati" can be parsimoniously used. To use this kind of market segmentation information, the analyst either has to impute the cluster membership of a customer from available information or has to estimate the probabilities of membership in each cluster and then apply Bayesian model averaging (Clyde and George, 2004) or some other ensemble method (Friedman and Popescu, 2005). Depending on the situation, the available information may not be sufficient to make a strong determination of cluster membership.

The first CRM task is to acquire a new customer: this is usually more expensive than retaining a current customer (five times more expensive, according to MBA folklore), and businesses want to target their recruitment investment to cherry-pick the most profitable ones. Businesses address this by combining statistical models of customer segments with individual information from cookies, purchased address lists and data warehouses.

The second CRM task is to please the customer; depending on the sector and the kind of customer, this may involve loss leaders, help-desk personnel, 24-hour service, information technology support and development for new web site services, and e-personalization (such as birthday cards). However, as the Wicked Witch of the West says, "These things must be done delicately or you hurt the spell." Some customers can be put off by overattentive service and alarmed that their on-line music provider knows their birthday; other customers are charmed.

The third CRM task is to retain customers. This is a moving target; competitors are constantly offering new services and prices. Data mining is a strong asset for the retention problem, over and above the usual advantage of customer inertia. Most businesses have detailed files on their customers that indicate all the user-account information, plus whatever other data have been acquired. Obviously, these records identify customers whose cell phone contract must be renewed or whose credit card is about to expire—these are the ones most likely to churn (switch) to a competitor. The records can also flag people whose spending habits are changing (a move, a new job, a nostalgic passion for comic book collecting). To the extent that data mining can interpret such changes, it allows the company to please customers by offering new services or payment plans designed for their current circumstances.

### 3.1 Trust

Customer trust is key to CRM. There are many ways to build it and even more ways to lose it. Businesses spend a lot of money on trust.

One issue is perception. Most (but not all) segments of consumers tend to give more trust to standard brands, to big companies that are not scandal-ridden and to those enterprises that seem large and rich. This is why companies have public relations managers, conscientiously polish an image and push quirky ad campaigns that ensure name recognition (such as the ING commercials that publicize the name, not the product).

A second issue is professionalism. If the web site works conveniently, if there is a logo and if the software does not glitch, then users develop confidence in the company's competence, and building this kind of across-the-board competence requires investment.

The third issue is the main one—security. A major hurdle to early e-commerce was giving people the confidence to make credit card purchases on-line. People still do not know the details behind encryption certificates and third party collection (e.g., PayPal), but they have generally learned how to take appropriate precautions. The recent spate of spam purporting to come from PayPal and e-Bay is an object lesson on how ephemeral a good reputation can be. It also shows that reputation is an exploitable asset; the scammers are faking PayPal's address precisely because its reputation allows them to gain and abuse access to financial data.

Security is broader than payment systems. Intrusion detection is essential for corporations: if their system is hacked and their data are misused, it may be impossible to recover public trust. Intrusion detection is a hard problem that has three versions—two statistical and one social. First, commercial methods such as McAfee VirusScan use pattern matching (signature detection) to identify known attacks; this protects unless the attack is deeply disguised. Second, anomaly detection is used to identify hitherto unseen attacks or to flag suspicious activities by an authorized user. Anomaly detection systems need a statistical model of normal transaction behavior: if the current traffic deviates too much (a decision that may depend on time series models with exponentially weighted moving averages to downweight the past), then the anomaly detector reports a possible attack (cf. Maxion and Tan, 2002). The third kind of security refers to everything from insider malfeasance to pharming and phishing (in pharming, the hacker exploits a vulnerability in the DNS (Domain Name Server) control software to capture a domain site and divert traffic from it to a false site, often with the intent of phishing, which entails tricking people into giving up confidential information). Although in some cases statistics may be helpful in discovering such violations, it does not usually play a large role.

Even more than viruses, Web customers fear inappropriate use of personal information; that is why many vendors list privacy policies on their web sites.

Similarly, people have different levels of comfort with cookies, which enable functionality that consumers value, but have not always been used with integrity. Also, the use of a commercial web site should not allow viruses or worms to access user systems.

Companies invest in building and publicizing secure systems for e-commerce. They may also adopt privacy policies that forswear a secondary income stream from marketing customer information (or offer discounts to customers who waive certain kinds of privacy restrictions).

Another way to build trust is to help people. This is one of the motives behind cooperative efforts by many e-commerce giants to "Can Spam!" (another motive might be to forestall the situation that PayPal and e-Bay have encountered). Before these highly publicized efforts, data miners were making large contributions. In the ongoing arms race, the best of today's spam filters are Bayesian. They use bag-of-words models to develop classification rules and claim very high success rates (cf. Madigan, 2005). Probably support vector machines and random forests could also be successful, but without such statistical technology, the Internet would be a more barren marketplace.

In that same spirit, we have intrusion detection, which also has statistical aspects. Commercial methods such as McAfee VirusScan use pattern matching (signature detection) to identify known attacks. Similarly, anomaly detection is needed to protect against unknown attacks. Research systems such as Harbinger (Maxion) rely on robust time series models and exponentially weighted moving averages of packet statistics to quickly find and flag suspicious transactions.

One apparently unresearched application of data mining is to understand the trade-offs in different trust-building strategies. Among the young digerati, an awkward web site costs business, but among back country folk, privacy may be the dominant concern. Companies can use statistical analyses to decide which kinds of improvements are cost-effective for their target clientele.

### 3.2 Efficiency

At a metalevel, businesses need to allocate their resources adaptively across all aspects of CRM. Investing everything in customer recruitment and nothing in customer satisfaction is a myopic path to bankruptcy. However, the problem of optimal investment is hard to solve. Notable aspects of this problem include the following:

- Product and web site development have large start-up costs, but once made, these can often be maintained for much less.
- Acquiring customer information has random payoffs. One does not know before the model is built whether it will successfully drive the ad campaign or the product line or the retention effort.
- Return on investment changes over time: for example, one runs out of cherry-picked customers, and the next layers are less profitable; new fads emerge; the cost of fuel or financing varies over time.
- Markets are adversarial. The front-runner is targeted by competitors, who may operate at a loss to woo new customers or collude to undermine the leader (especially in the realm of global e-commerce, where regulation is slender).
- Different strategies are needed at different phases of growth. The transition from a start-up with angel funding to a middle-sized company to an economic powerhouse is extraordinarily complex, and the people who can manage one phase well think differently from those who succeed in another.

The statistical tool to address such problems is portfolio analysis, but uncertainties and complexity overwhelm formal solution. The practical approach is probably game-based simulation, which would entail the development of commerce versions of red teaming (as used by the military for scenario exploration) or SimCity (a game based on resource trade offs in city planning).

Such simulation tools are not yet available, at least not in any formal, routinized way. Businesses therefore decide CRM investments more robustly, through rules of thumb, instinct and necessity. This creates inefficiencies in capitalization and misinvestment that can be exploited by competitors who are smarter—or luckier.

Good luck can be disastrous for a lightly capitalized start-up. A simplified version of the Barings Bank debacle shows why. Suppose Barings hired 1,000 new traders and suppose each adopted a risk strategy that gave them a 10% chance of being a star at the end of a year. About 100 traders were lucky: they were lauded for their insight, promoted and given authority to manage larger accounts. Hence they repeated their risk strategy: about 10 were successful in the second year; in the third year, Nick Leeson was the sole winner. Barings Bank empowered their *wunderkind* to make large investments in derivatives futures; his next gamble, with statistical inevitability, destroyed them. The point of the parable is that it is easy for people to conflate luck with skill. Wise executives need to shrink toward the mean when assessing their own success, and the managers of data miners need to keep this in mind as well, since many discoveries will be spurious.

## 4. BUSINESS STRATEGIES FOR THE FUTURE

Many business practices can benefit from mining e-commerce information, even if they are not directly using it to promote new services or better handle their customers. Information collected by e-commerce transactions can inform businesses that practice no e-commerce at all.

An immediate example is the use of e-Bay auction data to determine the price points for certain kinds of products. That knowledge can affect both the marketing strategy and the mix of products that are manufactured. Shmueli and Jank (2005) analyzed this kind of information and even provided graphical summaries of the bidding dynamics. There is further work to be done; eBay now allows the seller to set a price which, if met, terminates the auction. Analyzing the exercise of that early ending option gives insight to the distribution of utility values in the e-Bay population.

As another example, companies that do international business need constant guidance on ever changing regulations, tariffs, price differentiation across geography, weather conditions, lading fees and spot prices for commodities. This field is complex and highly localized; one needs smart people on the ground in every country with which business is done. It is expensive to build this kind of expertise in-house, and probably impossible for it to be truly expert, since complex optimization problems must be solved. As an e-business opportunity, a start-up could hire locals to input this kind of information from all over the world into a common data base. Embellished with appropriate data mining tools, the company could then sell access to this planning system to all international shippers.

Data mining could create a sales force that seemed astonishingly impressive. Imagine that an investor called a broker about buying stock in XYZ, Inc. If the broker promptly replied that XYZ could be an attractive investment because it was the only local fertilizer company in Connecticut that had shown a profit in more than half of the last six years, the investor would be dazzled. Of course, the broker's trick is to have a search system that can find *something* worthy of comment about almost any stock—the curse of dimensionality ensures that almost every case is an extreme point in high dimensions. Such specialized search systems seem straightforward to build and would have applications beyond generating a facade of business acumen (which smart repeat buyers would eventually learn to discount). Sports enthusiasts and other trivia buffs would be an avid audience.

Besides these exotic possibilities, there are more standard applications:

- Transaction mining can indicate which stores need which product (Victoria's Secret uses this to send specific apparel to specific stores).
- One can create virtual tours of real brick-and-mortar stores to help entice shoppers and support more human-oriented searches (Broadvision.com helps to set up such web sites; Powell's Books used to show shelf photos, but it appears to be aping Amazon now).
- Association rules can determine which items are purchased together, leading to strategic layouts of stores or artful discount strategies; for example, drop the price on spaghetti sauce, but make it back on pasta. (The tiredest example concerns beer and diapers; all that is going on here is that houses with babies tend to have a resident male between 20 and 40 who is not going out to the bars very much.)
- Recommender systems can now support books, music and movies; some cluster analysts believe they can even identify the characteristic features of a hit pop tune.

There are many other examples. The point is that the e-commerce toolkit has broad application and can play a role in applications far beyond the stereotype of a customer sitting with a laptop.

The potential for linking lots of hard data with really sophisticated algorithms is enormous. Only a few companies now are truly pushing the envelope of high-end statistical analysis, but as the level of performance rises, everyone will have easier access to flexible tools and experts with the knowledge to use them. This will enable businesses to become presbyopic: they can plan better, control better and adapt faster. There is reason to hope that the well-documented inefficiencies of the market will be reduced, which economists suggest will make the world better (in the long run).

Market efficiency may be what businesses should fear most. Suppose customers had access to the kind of high-end data mining tools that e-commerce uses. People could use intelligent personalized shopping bots to churn endlessly, always demanding the best value. This would drive down profit margins by forcing tight competition. For example, the bot could adaptively review cell phone use history before each call, forecast costs and switch automatically to the plan that is most advantageous. Businesses could discourage this by offering long-term plans at lower rates, but this creates market niches for microplans or low-fee long-term agreements.

It is easily conceivable that electronic commerce will someday lead to micro-haggling on every transaction; Priceline.com is already moving in this direction. At the end of the day, we may see a business community in which all margins are razor thin and every transaction has a specific cost.

The analysis of such questions is important. Electronic commerce provides an unusual setting for statisticians and decision-makers: it is rich in data, but poor in plausible assumptions that are useful for modeling, such as independent observations. Our challenge is to use the abundance of data in place of tractable simplicity, but still extract interpretable descriptions that support business decisions. This is made more difficult by the expansion of e-commerce, because more people are shopping on-line and more businesses are adding electronic services, so any analysis must track a moving target.